\documentclass{article}
\usepackage{graphicx}
\usepackage[numbers,sort&compress]{natbib}
\usepackage{amsmath,amsthm}
\usepackage{authblk}

\def\Re{{\rm Re}}

\title{The approximate conformal mapping of a disk onto domain with an acute 
angle }

\author[1]{P.~N.~Ivanshin}
\author[1]{E.~A.~Shirokova}
\affil[1]{Lobachevskiy Institute of Mathematics \& Mechanics, Kazan Federal
University, Kremlevskaya st., 35, Kazan, 420008, Russia.
E-mail: pivanshi@yandex.ru}

\begin{document}
\maketitle

\begin{abstract}

The  method of boundary curve reparametrization is applied to construction of the approximate analytical conformal mapping of the unit disk  onto an arbitrary given finite domain with a boundary smooth at every point  but fininte number of  acute angle points. The method is based on both the Fredholm equation solution and spline-interpolation. This approach consists of approximate solution of a linear system with unknown Fourier coefficients and construction of correction  splines. The approximate  mapping function has the form of a Cauchy integral. The method presentation is supported by demonstration of some examples. This method is applicable to the case of multiply connected domains with boundary angle points.


\end{abstract}

\section*{1. Introduction}

Conformal mappings by the analytical functions
of complex variable play an important role in solution of many problems
of mechanics and mathematics, particulary in the case of plane potential
fields and Laplace equation solution \cite{1}. The conformal mapping of
the circular domain onto a given domain with complex boundary can be
applied to solution of plane boundary value problems for corresponding
domains by symmetry methods, for example, with the help of Schottky
-- Klein prime functions known for a
circular domain \cite{2}. Computer progress stimulated appearance of many
numerical methods for conformal mapping construction \cite{3}. All these
methods are rather time-consuming. For example, the widely used Wegmann
numerical method is based on Riemann--Hilbert problem solution
and involves iteration processes \cite{4, 5}.

There are several types of canonical regions for conformal mappings
\cite{6}. The five types of canonical slit regions are disk with concentric circular
slits, annulus with concentric circular slits, unbounded circular slit regions,
unbounded radial slit regions, and unbounded parallel slit regions. Nasser
presented the method of mapping construction of bounded and unbounded multiply connected regions
onto these five canonical regions by reformulating the mapping function
as a Riemann-Hilbert problem which is solved by means of boundary
integral equation with the generalized Neumann kernel. The right-hand
side of the integral equation involves integral with cotangent singularity
which is approximated by Wittich's method. The integral equation was
discretized by the Nystr\"{o}m method with the trapezoidal rule to obtain a
linear system \cite{7, 8, 9, 10, 11, 12, 13}. The method presented in \cite{6,7,8,14,15,16} is based also on the approximate solution of the integral equation. There the mapping construction is reduced to solution of a truncated linear system over the unknown Fourier coefficients and the problems with the cotangent singulariy are overpassed with the help of the Hilbert formula. 
 
Some of the problems of physics are connected with domains with boundary angle points. It is not difficult to find an 
appropriate mapping onto such domain by application of the integral equation method  when the inner angles of the domain are greater than $\pi$. At the same time,  for the case of the  boundary angle  point with the inner angle less than $\pi$ the kernel of the corresponding  Neumann equation has the singularity  at  this point \cite{Pm2, Iar}. The numerical solution of the mapping problem for these domains  is also difficult. Some work is done regarding the optimal mesh structure of   $n$-gonal domains \cite{BB}.   Certain authors considered asymptotic expansions with logarythmic terms \cite{Leh}. D. M. Hough and N. Papamichael  dealt with polygonal domains gluing splines and singular functions \cite{Pm2}.  Construction of conformal mapping for  domains with angle points via application of the additional conformal mappings and the continued fractions  was presented in \cite{Iar}. 

%At the same time, domains with inner angles less than $\pi$ and especially, acute, present much greater difficulty. Note also that the case of acute boundary angles is closely related to the case of so-called thin domains.    Note that the only way to work in the neighbourhood of the angle point in these conditions is a dense node distribution%

Here we present a  method of the approximate conformal mapping
of the unit disk 
onto a  simply connected domain with acute angle points at its boundary applying both trigonometric polynomials and splines.  

We first construct the mapping of the unit disk onto a similar domain with the  smooth  boundary. We apply the standard method based on solution of the integral equation over the imaginary part of the boundary value of an auxiliary regular function  and the linear system truncation.Then the analytic function is restored via Cachey integral. This construction is given in Section 2. The method is adequate for the case of a smooth boundary of the domain. Presence of an acute angle boundary point of the given domain brings to the linear system solution disturbance. Therefore the dependence of the polar angle of the  unit circle on the boundary parameter of the given domain fails to be monotonic.

After we construct the mapping with nonmonotonic polar angle of the unit circle we consider the piecewise correction of this polar angle function with the help of monotonic splines in the neighbourhoods of the preimages of the angle points. Now the Cauchey integral is constructed with application of  reparametrization of the boundary parameter with the inverse of spline function.  Here we have no need in gluing together two sets of solutions and simply modify the already constructed one as in \cite{Pm2}.  This is the content of Section 3.

\section*{2. Approximate conformal mapping of the unit disk  onto a  domain similar to the given one by means of integral equation}

Consider a simply connected domain $D_z$ bounded by the simple
 curve $ L_0$  given by the equation $z=z_0(t)$, $t\in [0,2\pi]$.

We also assume that the boundary curve $L_0$ complex representation is as
follows:
$$
z_0(t)= \sum\limits_{k=-m}^{n}d_{k } e^{i k t},\ t\in [0,2\pi].
$$
The parametrization traces the domain $D_z$ along $L_0$ counterclockwise.

%
%Definition 1. We call the unit  disk with $m$  circular slits $\zeta = R_s e^{i \theta}$, $\theta_{1 s} <\theta <\theta_{2 s}$,  $0<R_s<1$, $\theta_{2 s} - \theta_{1 s} <2\pi$, $s=1,\ldots,m$,  and with $n-m$ radial slits $\zeta = R e^{i \theta_j}$, $0<R_{1 j} <R <R_{2 j}<1$,  $0<\theta_j<2\pi$,  $j=m+1,\ldots,n$, an $(n+1)-$  connected canonical domain of the first type.
%
%Definition 2. We call the annulus with the exterior radius 1,  with  the interior radius $r$, $r<1$, and with  $(m-1)$  circular slits $\zeta = R_s e^{i \theta}$, $\theta_{1 s} <\theta <\theta_{2s}$,  $r<R_s<1$, $\theta_{2 s} - \theta_{1 s} <2\pi$, $s=1,\ldots,m-1$, and with $n-m$ radial slits $\zeta = R e^{i \theta_j}$, $0<R_{1 j} <R <R_{2 j}<1$,  $0<\theta_j<2\pi$,  $j=m+1,\ldots,n$, an $(n+1)-$  connected canonical domain of  the second type.
 We present  construction of the mapping of the unit disk  onto  a domain similar to  $D_z$ due to \cite{7} . 

%\begin{theorem}  
%
%There exists  a regular function $z=f(\zeta)$ in the unit disk $D$ such that the function $f(\zeta)$ maps  conformally  the domain $D$ onto the given  domain $D_{z}$  with smooth boundary. 
%%the image of the right horizontal slit boundary point is the point $z_0(0) \in L_0$. 
%\end{theorem}
%
%
%\begin{proof}   

We assume that $0\in D_z$  without loss of generality. 
  %We construct the conformal map of the unit disk onto the domain $D_z$ by reparametrization of the given boundary representations. So we search for the function $t(\theta),$ $\theta \in [0, 2\pi]$.

Let us  consider the analytic in the domain $D_z$ function $\zeta(z)$ which maps conformally the
domain $D_z$ onto the unit disk $D_{\zeta}$ with the correspondence $\zeta(0) = 0$  and the 
analytic in $D_z$ function $\psi(z)=\log \frac{z}{\zeta(z)}$. Since the necessary and sufficient
condition for  $\psi(z_0(t))=\log \frac{z_0(t)}{e^{i\theta_0(t)}}$ to be the boundary values of the analytic in $D_z$ function we have the boundary relation
$$
\log \frac{z_0(t)}{e^{i\theta_0(t)}} = \frac{1}{\pi i} \int\limits_0^{2\pi}
\log \frac{z_{0}(\tau)}{e^{i\theta_0(\tau)}} [\log(z_{0}(\tau)-z_0(t)) ]'_\tau d\tau.  \eqno(1)
$$
where $ t\in[0, 2\pi]$.

%We introduce the following functions: $q_s(t) =\arg z_s(t) - \theta_s(t)$,  where $\theta_s(t)$ is the polar angle of the image of the point of $z_s(t)$, $s =0, \ldots, m$, and $p_j(t)=\log |z_j(t)| -\log R_j(t)$, where $R_j(t)$ is the radius of the image of the point of $z_j(t)$, $j =m+1, \ldots, n$. 

We separate the imaginary part of both sides of equation (1):
$$
q_0(t)= \frac{1}{\pi} \int\limits_{0}^{2\pi} q_{0}(\tau) [\arg(z_0(\tau)-z_0(t)) ]'_\tau  d\tau -
$$
$$
 - \frac{1}{\pi}
\int\limits_{0}^{2\pi} \log |z_{0}(\tau)| [\log|z_0(\tau)-z_0(t)| ]'_\tau  d\tau, \eqno(2)
$$
where $q_0(t)=\arg(z_0(t))-\theta_0(t)$, $\theta_0(t)$ is the dependence of the polar angle of the unit circle on the boundary parameter of the domain $D_z$.

%$$
%+\sum\limits_{j=m+1}^n \frac{1}{\pi} \int\limits_{0}^{2\pi} [\arg z_j(\tau)-\theta_j] [\arg(z_j(\tau)-z_s(t)) ]'_\tau  d\tau -
%$$
%$$
% - \sum\limits_{j=m+1}^n \frac{1}{\pi}
%\int\limits_{0}^{2\pi} p_j(\tau) [\log|z_j(\tau)-z_s(t)| ]'_\tau  d\tau,  \quad s=0, \ldots, m.
%$$

After differentiating  relation (2)  with respect to $t$ and integrating the results by parts, we
obtain the following relation on the function $q_0'(t)$:
$$
q_0'(t)=\frac{1}{\pi} \int\limits_{0}^{2\pi} q_{0}'(\tau) K(\tau,t)  d\tau+ P(t),\eqno(3)
$$
%$$
%+\sum\limits_{j=m+1}^n\frac{1}{\pi} \int\limits_{0}^{2\pi} p_j'(\tau) L_{j, s}(\tau,t)  d\tau  \ \ s=0, \ldots, m, 
%$$
where
$$
 K(\tau,t) = - [\arg(z_0(\tau)-z_0(t)) ]'_t,\ \ L(\tau,t)=[\log|z_0(\tau)-z_0(t)| ]'_t,
 $$
% $$ 
%Q(t)=  \frac{1}{\pi}
%\int\limits_{0}^{2\pi} [\log |z_{0}(\tau)|]' L(\tau,t) d\tau, 
%$$
$$ 
P(t)=   \frac{1}{\pi}
\int\limits_{0}^{2\pi} [\log |z_{0}(\tau)|]'_\tau L(\tau,t) d\tau. 
$$

The kernel $L$  has a singularity in the form of $\cot \frac{\tau-t}{2}$:
$$
(\log|z_0(\tau)-z_0(t)| )'_t =\Re \left(\log  \sum\limits_{k=1}^{n}d_{k }[e^{i k\tau}-e^{i kt}]\right)'_t=\Re\left(\log \sin\frac{\tau-t}{2}+\right.
$$
$$
\left.+\log\left[ \sum\limits_{k=1}^{n} d_{k}e^{ikt}  \sum\limits_{l=0}^{k-1} e^{i l (\tau-t)}- \sum\limits_{k=1}^{m_s} d_{(-k) }e^{-ik\tau}  \sum\limits_{l=0}^{k-1} e^{i l (\tau-t)}\right]\right)_t'=
$$
$$
=-\frac{1}{2}\cot \frac{\tau - t}{2} + \left(\log\left|\sum\limits_{k=1}^{n} d_{k }e^{ikt}  \sum\limits_{l=0}^{k-1} e^{i l (\tau-t)}- \sum\limits_{k=1}^{m} d_{(-k) }e^{-ik\tau}  \sum\limits_{l=0}^{k-1} e^{i l (\tau-t)}\right|\right)_t'.
$$

The Cauchy principal value integral
$$
\frac{1}{\pi}
\int\limits_{0}^{2\pi} [\log |z_{0}(\tau)|]' \cot\frac{\tau-t}{2} d\tau
$$
can be calculated via Hilbert formula \cite{9} as in \cite{8}.

If we search for the  solution of equation (3) in the form of Fourier series:
$$
q_0'(t)= \sum\limits_{l=1}^{\infty} \alpha_{l } \cos lt+ \beta_{l } \sin lt, \
$$
given
$$
P(t)= \sum\limits_{l=1}^{\infty} \gamma_{l } \cos lt+ \kappa_{l } \sin lt, \  \ t\in[0,2\pi],  
$$
we can rewrite equation (3) as  the following infinite system of  equations  which can be written in the operator form as follows:
$$
\left(\begin{array}{cc}
I- K_{c,c} & -K_{c,s}\\
- K_{s,c} & I- K_{s,s}  
\end{array}\right) \left(\begin{array}{c}\alpha \\
\beta \\
\end{array}\right) 
=\left(\begin{array}{c}
\gamma \\
\kappa \\
\end{array}\right). \eqno(4)
$$
Here $$
\alpha=(\alpha_1,\alpha_2,\hdots)^T,\      \beta=(\beta_1,\beta_2,\hdots)^T,\     \gamma=(\gamma_1,\gamma_2,\hdots)^T,\    \kappa=(\kappa_1,\kappa_2,\hdots)^T  .
$$

The last operator system is in fact the infinite linear system over the Fourier coefficients of the unknown function $q_0'(t)$,  if we find the coefficients of double Fourier expansions of  the kernels of integral operators and compare the coefficients with the same trigonometric functions  \cite{6}. 

The solution $q'_0(t)$ of system (4) allows to obtain the polar angle function: $ \theta_0(t)=\arg(z_0(t))-q_0(t)$. Note that $\theta_0(t)$ should grow monotonically when $t$ grows from $0$ to $2\pi$, $\theta_0(2\pi)-\theta_0(0)=2\pi$. Now the Cauchey integral provides the analytic function $f(\zeta)$ that maps the unit disk onto the given domain $D_z$:
$$
f(\zeta)=\frac{1}{2\pi} \int\limits_{0}^{2\pi}\frac{z_0(t) e^{i\theta_0(t)} \theta_0'(t) d t}{e^{i\theta_0(t)}-\zeta}.\eqno(5)
$$
   
     Approximate solution of the infinite system over Fourier coefficients of the unknown functions is a solution of a truncated system over the Fourier coefficients of the unknown functions.

We search for the approximate solution of system (4) in the form of a Fourier polynomial:
$$
\tilde{q}_0'(t)= \sum\limits_{l=1}^{M} \alpha_{l } \cos lt+ \beta_{l } \sin lt, \  \ \eqno(6)
$$
given
$$
\tilde{P}(t)= \sum\limits_{l=1}^{M} \gamma_{l } \cos lt+ \kappa_{l } \sin lt, \  \ t\in[0,2\pi].  
$$

Now integral Fredholm equations of the second kind in (3) can be reduced to
the finite linear system over Fourier coefficients $\alpha_{j}$ and $\beta_{j}$:

$$
\left( \begin{matrix}A & B \\
C & D \end{matrix} \right) \cdot
 \left( \begin{matrix} \alpha_1 \\  \vdots\\ \alpha_M \\ \beta_1\\ \vdots\\ \beta_M \end{matrix} \right) =
\left( \begin{matrix} \gamma_1 \\ \vdots\\  \gamma_M \\ \kappa_1\\ \vdots \\ \kappa_M \end{matrix} \right).  \eqno(7)
$$
The vectors   $(\gamma_{1 },\ldots,\gamma_{M })^T$, $(\kappa_{1 },\ldots,\kappa_{M})^T$ on the right-hand side of the system consist of the elements
$$
\gamma_{j }=\frac{1}{\pi}
\int\limits_{0}^{2\pi} P(t)\cos jt dt,\  \kappa_{j }=\frac{1}{\pi}
\int\limits_{0}^{2\pi} P(t)\sin jt dt,\  \ j=1, \ldots, M.
$$

The block matrices $A$, $B$, $C$, $D$, of size $ M\times M$ consist of the elements
$$
A_{j k} =  \delta_{j k} - \frac{1}{\pi^2}\int\limits_{0}^{2\pi}\cos k\tau d\tau  \int\limits_{0}^{2\pi}K(\tau, t) \cos jt dt, 
$$
$$
B_{j k} =  - \frac{1}{\pi^2}\int\limits_{0}^{2\pi}\sin k\tau d\tau  \int\limits_{0}^{2\pi}K(\tau, t) \cos jt dt, 
$$
$$
C_{ j k} =  - \frac{1}{\pi^2}\int\limits_{0}^{2\pi}\cos k\tau d\tau  \int\limits_{0}^{2\pi}K(\tau, t) \sin jt dt, 
$$
$$
D_{j k} = \delta_{j k} - \frac{1}{\pi^2}\int\limits_{0}^{2\pi}\sin k\tau d\tau  \int\limits_{0}^{2\pi}K(\tau, t) \sin jt dt,
$$
where $ j, k=1, \ldots, M$, $\delta_{j k}$  is the Kronecker delta function.

The approximation $\tilde{q}(t)$ of the function $q_0(t)$, can be restored via the derivative (6) with
an arbitrary constant summand:
$$
 \tilde{q}(t)=\sum\limits_{l=1}^{M} \frac{\alpha_{l }}{l} \sin lt- \frac{\beta_{l }}{l} \cos lt, 
$$

%So the parameter of the  domain $D_{\zeta}$  has been found.

Now we have the approximation $\tilde{q}(t)$ of the function $q_0(t)$,  $t\in [0,2\pi]$, and therefore we can restore the approximate relation between the boundary parameter of the domain $D_z$ and the polar angle of the boundary of $D_{\zeta}$ via the formula $\tilde{\theta}(t)=\arg z_0(t)-\tilde{q}(t)$.  

The approximate analytical function which maps $D_{\zeta}$ onto $D_z$ now has the form of the Cauchy integral by analogy with (5):
$$
\tilde{f}(\zeta)=\frac{1}{2\pi} \int\limits_{0}^{2\pi}\frac{z_0(t) e^{i\tilde{\theta}(t)} \tilde{\theta}'(t) d t}{e^{i\tilde{\theta}(t)}-\zeta}.\eqno(8)
$$

%We can apply the Cauchy integral in the form 
%$$
%f(\zeta)=\frac{1}{2\pi} \int\limits_{0}^{2\pi}\frac{z_0(t) e^{i\theta_0(t)} \theta_0'(t) d t}{e^{i\theta_0(t)}-\zeta}
%$$
%in order not to deal with the functions $t_s^{\pm}(\theta)$ or $t_j^{\pm}(R)$ and not to integrate alond the different borders of the same slit.

%\end{proof}

Existence of the exact solution of system (4)  and convergence of the approximate
solution of system (7) to the exact one provided $M \to \infty$ were proved in \cite{7} for the case of the smooth boundary of   $D_z$.

Since we  reduce integral equation system (3)   first   to infinite linear system  (4) and then to the finite one given by relation (7) we need to describe the behaviour of the coefficients of the infinite and the truncated  matrices of the linear systems (4) and (7) in the case when  the boundary of $D_z$ posesses an  angular point. The elements of the infinite matrices $K_{cc}$, $K_{cs}$, $K_{sc}$ and $K_{ss}$ in (4) are the double Fourier series coefficients of the kernel   $K(\tau,t)=-[\arg(z_0(\tau)-z_0(t)) ]'_t$. This kernel is not continuous  at the point $t_0$ if the boundary of $D_z$ contains the angle point $z(t_0)$ with the inner angle $\lambda \pi$ when $0< \lambda <1$ but has the  singularity of type $\frac{1}{|\tau-t|^{1-\lambda}}$,  \cite{Iar, Pm2}. This singularity disturbes the convergence of the  finite system (7)
%$(I_{2M}-P_{2M})X_1=Y_1$ 
to the infinite system (4) so that the corresponding    approximation  $\tilde{\theta}(t)$ of the polar angle $\theta_0(t)$ fails to be monotonic and therefore the function (8) can not provide an adequate mapping.

\section*{3. Reparametrization in the neighbourhood of the singular point}

We suppose that the boundary of $D_z$ contains the acute angle point at $z_0(t_0)$ with the inner angle $\lambda \pi$, $\lambda <1$. The first step in our construction of the approximate mapping function of the unit disk onto $D_z$  is approximate solution of integral equation (3). %actually the mapping from the unit disk onto the angled domain. 
%It is known \cite[Statement 1]{Iar} that then the kernel $K(\tau, t)$ of integral equation (3) possesses a singularity.

The approximation  $\tilde{\theta}(t) $ to $\theta_0(t)$ fails to be monotonic in the neighbourhood of the parameter value $t_0$, so  
 the normals to the  boundary of $D_z$ overlap at the neighbourhood of the angle point as it is demonstrated  %on the left side of  So, the reparametrization $\theta_0(t)$ generally has the form given%  
  by Fig. 1, a. In order to save the function monotonicity we cut the fold (see Fig. 1, b) in the following way. 

    \begin{figure}[h]\label{fig1}
   \begin{center}
  \includegraphics[width=4truecm,height=4truecm]{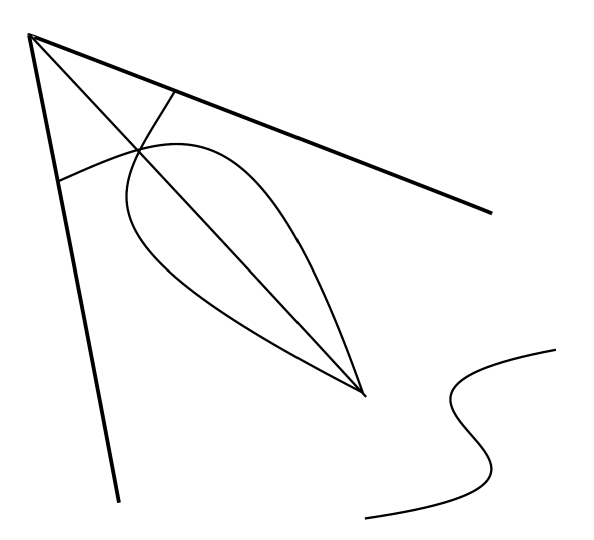}   
   \includegraphics[width=4truecm,height=4truecm]{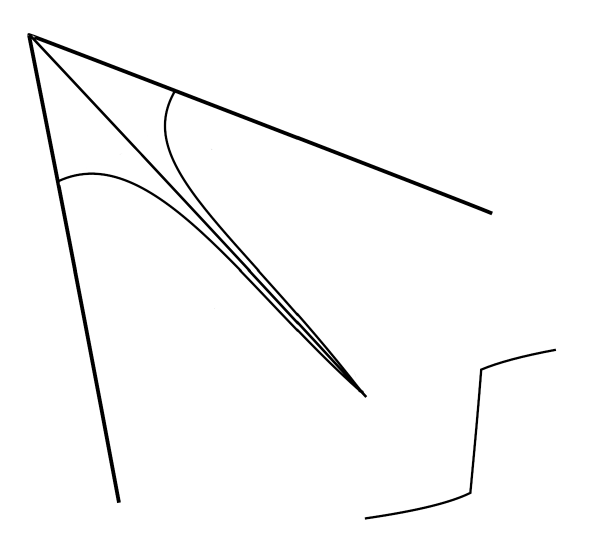} 
\caption{Normals overlap in the neighbourhood of the angle point}
  \end{center}
  \end{figure}

Note that the polalar angle $\theta_0(t)$ vanishes at the point $t_0$. Really, due to the presence of the angle point the mapping function in the neighbourhood of its preimage has the following representation: 
$$
f(\zeta)=f(e^{i\theta_0(t_0)}) +  (\zeta-e^{i\theta_0(t_0)})^{\lambda} g(\zeta),
$$
so 
$$
f'(\zeta)=\frac{\Psi(\zeta)}{(\zeta-e^{i\theta_0(t_0)})^{1-\lambda}}.
$$

 Since 
$$
f'(e^{i\theta_0(t)})= \frac{z'_0(t)}{i e^{i\theta_0(t)} \theta'_0(t)}
$$ 
the singularity at the parameter value $t=t_0$ appears only if $\theta'_0(t_0)=0$. 

Therefore we replace the  function $\tilde{\theta}(t)$ by a spline $\phi(t)$ on a segment $[t_0-\epsilon_1,t_0+\epsilon_2]$ so that $\phi'(t)>0$, $t\in [t_0-\epsilon_1,t_0) \cup  (t_0,t_0+\epsilon_2]$, $\phi'(t_0)=0$ and the continuous function 
$$
\breve{\theta}(t)= \left\{ \begin{array}{cc}  \tilde{\theta}(t),&  t \in [0,2\pi]\setminus [t_0-\epsilon_1,t_0+\epsilon_2],\\
\phi(t),&  t\in [t_0-\epsilon_1,t_0+\epsilon_2],
\end{array}\right.
$$
monotone increases on $[0,2\pi]$ (see Fig. 2).

  \begin{figure}[h]\label{fig2}
   \begin{center}
  \includegraphics[width=5truecm,height=3truecm]{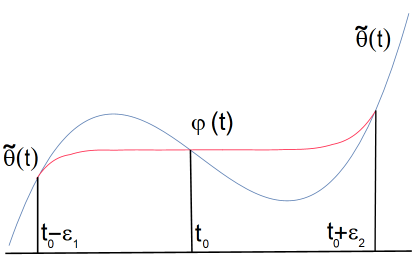}   
\caption{Spline in the neighbourhood of the angle point}
  \end{center}
  \end{figure}

Due to monotonic increase of the spline $\phi(t)$ we can construct the inverse also monotone increasing function $t=t(\phi)$, $\phi \in [\tilde{\theta}(t_0-\epsilon_1),\tilde{\theta}(t_0+\epsilon_2)]$. Now we  rewrite  formula (8) in the case of the acute angle point $t_0$:

$$
\hat{f}(\zeta)=\frac{1}{2\pi} \int\limits_{t_0+\epsilon_2}^{2\pi+t_0-\epsilon_1}\frac{z_0(t) e^{i\tilde{\theta}(t)} \tilde{\theta}'(t) d t}{e^{i\tilde{\theta}(t)}-\zeta}+ 
$$
$$
+\frac{1}{2\pi} \int\limits_{\tilde{\theta}(t_0-\epsilon_1)}^{\tilde{\theta}(t_0+\epsilon_2)}\frac{z_0(t(\phi)) e^{i\phi}  d \phi}{e^{i\phi}-\zeta}.      \eqno(9)
$$

Example 1. Unit disk with two orthogonal tangent lines. The left part of Fig. 3 is the result of the first step, cubic spline. The second step also with the cubic spline is the right part of Fig.3, red lines show the target domain boundary. The picture shows us that the approximation is better than in \cite{Iar}, the order of approximation being $O(0.001)$ in contrast with $O(0.01)$ of \cite{Iar} with the same size of the system (6).

  \begin{figure}[h]\label{fig3}
   \begin{center}
  \includegraphics[width=4truecm,height=4truecm]{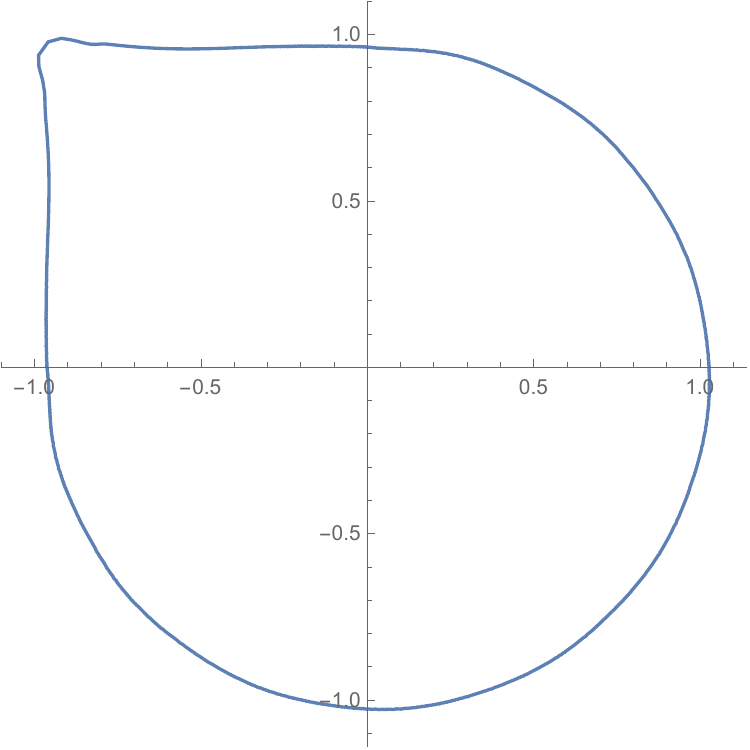}   
   \includegraphics[width=4truecm,height=4truecm]{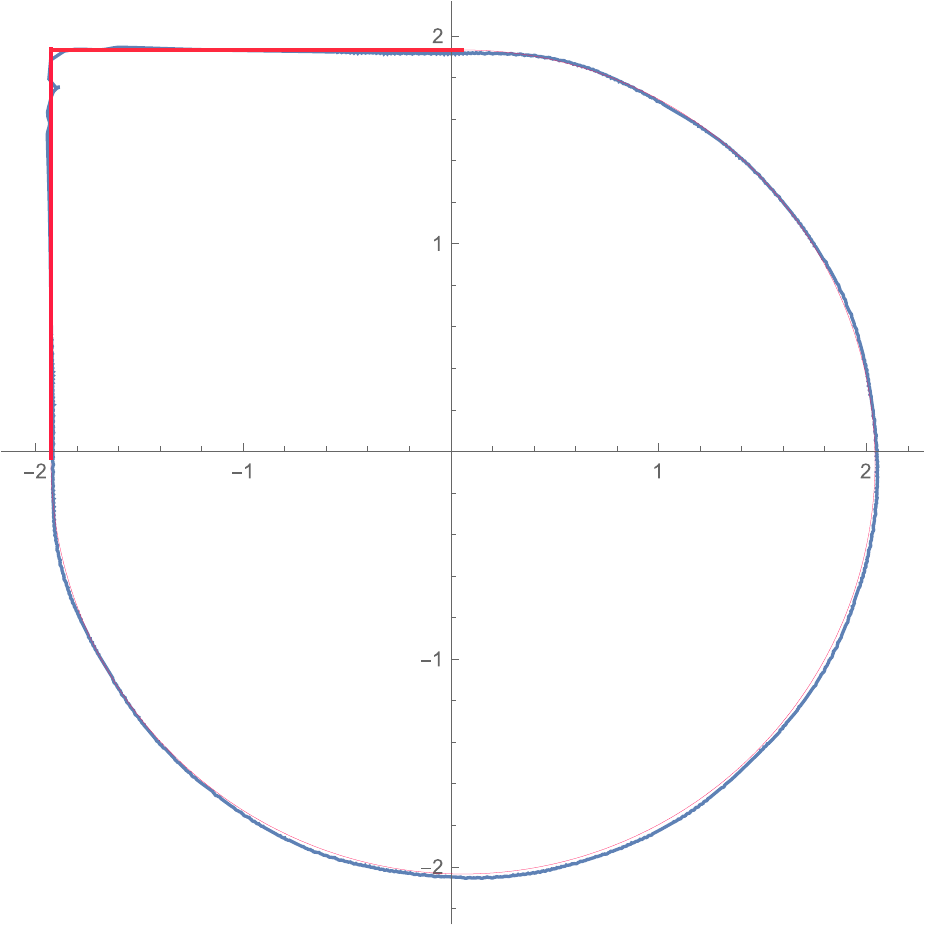} 
\caption{Domain with one angle}
  \end{center}
  \end{figure}

Example 2. A semidisk of radius $1$. The left part of Fig. 4 is the result of the first step, cubic spline. The second step with the linear spline is the right part of Fig. 3, red lines show the target domain boundary. The lower image is the angle reparametrization.

  \begin{figure}[h]\label{fig4}
   \begin{center}
 
   \includegraphics[width=4truecm,height=2truecm]{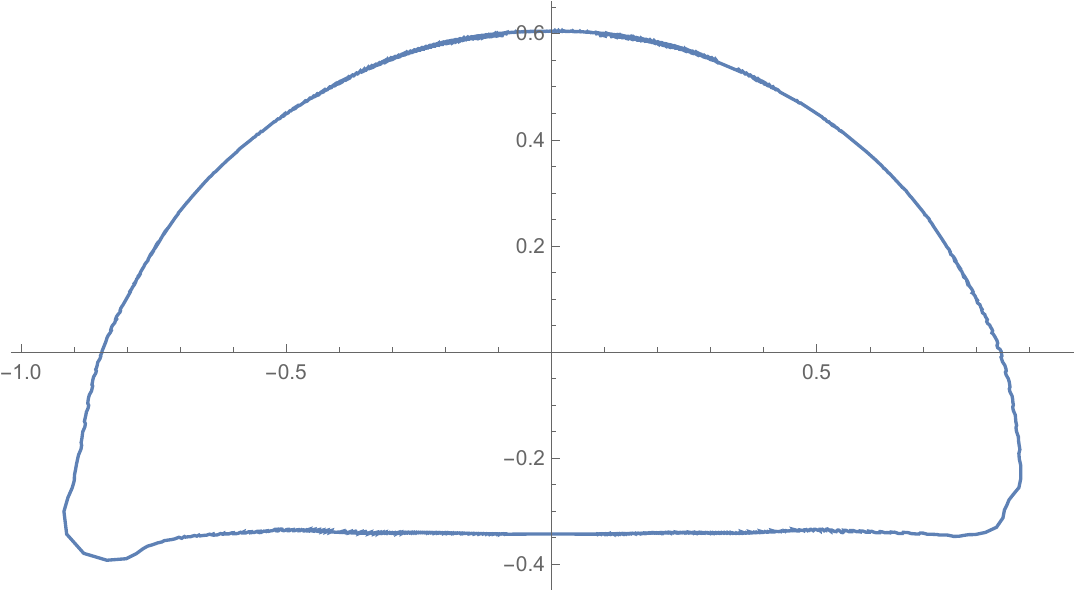} 
     \includegraphics[width=4truecm,height=2truecm]{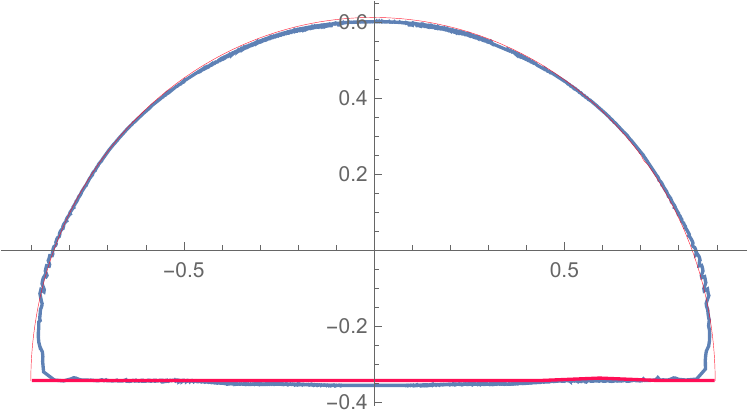}  
    \includegraphics[width=4.8truecm,height=4truecm]{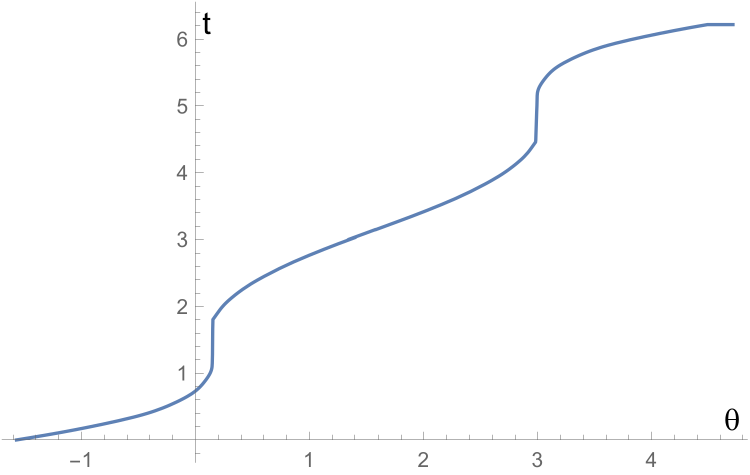} 
\caption{Semidisk and the angle reparametrization}
  \end{center}
  \end{figure}

Example 3. A figure bounded  by the curve $e^{i t}+\frac{1}{4} e^{-i 2 t}+\frac{1}{8} i e^{-i 3 t}$, $0 \leq t \leq 2 \pi$. It can be seen here that not only angles but also domains with relatively thin boundary elements can be approximated by the proposed approach. The first two parts of Fig. 5 are the results of the first step, cubic spline and the level line of the second step, linear spline; red lines show the target domain boundary. The lower image is the angle reparametrization.

  \begin{figure}[h]\label{fig5}
   \begin{center}    
   
   \includegraphics[width=3truecm,height=3truecm]{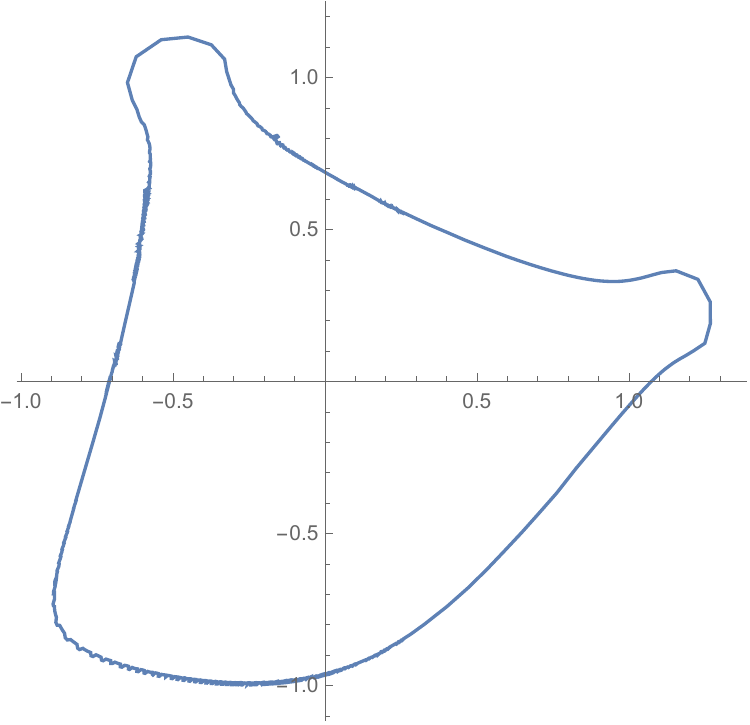}    
   \includegraphics[width=3truecm,height=3truecm]{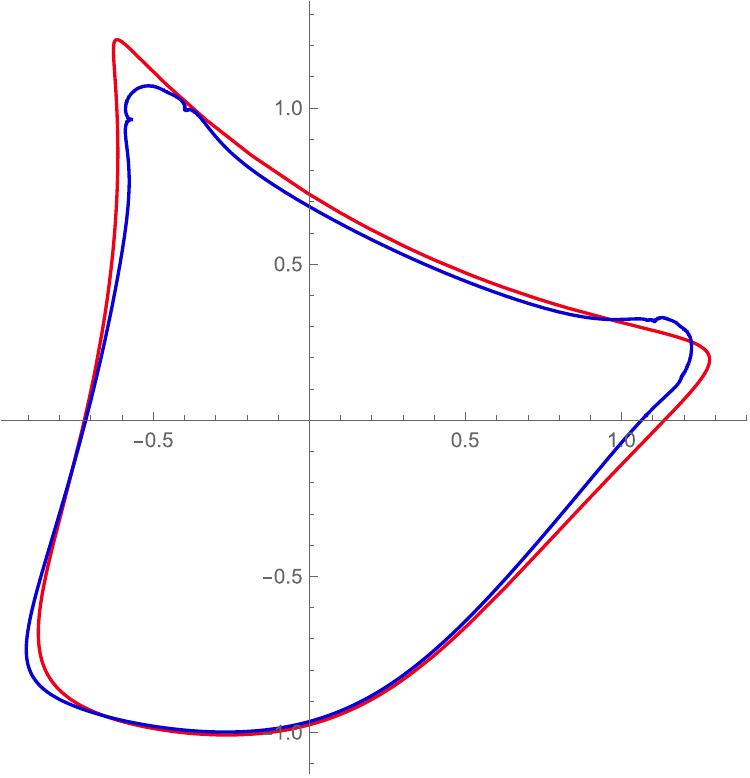}   \\
   \includegraphics[width=3.5truecm,height=3truecm]{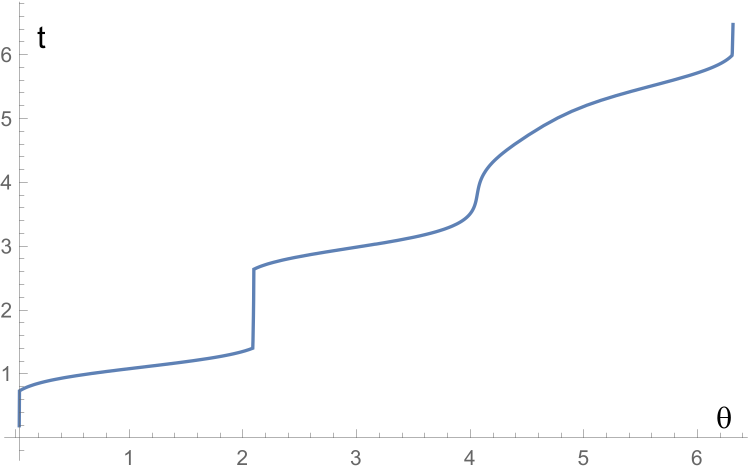} 
\caption{Two approximations and the angle reparametrization}
  \end{center}
  \end{figure}

\section*{Declarations}
 
{\bf Funding} not applicable.

{\bf Conflicts of interest}  We hereby state that the article ``The approximate conformal mapping of a disk onto domain with an acute 
angle'' does not involve any conflict of interests.

{\bf Availability of data and material} On demand.

{\bf Code availability} Fortran code and .nb available on demand

{\bf Authors' contributions} Both authors 50 per cent each.


\begin{thebibliography}{99}

\bibitem{1} R. Schinzinger, P. A. A. Laura. Conformal Mapping: Methods and Applications. Dover Publications, 2003.

\bibitem{2}  D. Crowdy, J. Marshall. Conformal Mappings between Canonical Multiply Connected Domains. Computational Methods and Function Theory, V.6, pp. 59-76, (2006).

\bibitem{3} A. H. M. Murid, Laey-Nee Hu. Numerical conformal mapping of bounded multiply connected regions by an integral equation method. Int. J. Contemp. Math. Sci., V.4, pp. 1121-1147 (2009).

\bibitem{4} R. Wegmann. Fast conformal mapping of multiply connected regions. J.Comput. Appl. Math., V.130, pp. 119--138 (2001).

\bibitem{5} R. Wegman, M. M. S. Nasser.  The Riemann-Hilbert problem and the  generalized Neumann kernel on multiply connected regions, J. Comput. Appl.  Math., V.214, pp. 36--57 (2008).


\bibitem{6} E. A. Shirokova. On the approximate conformal mapping of the     unit disk on a simply connected domain. Russian Math., V.58, N.3,   pp.  47--56, (2014).

\bibitem{7} E. A. Shirokova, P. N. Ivanshin. Approximate Conformal Mappings and
Elasticity Theory. J. of Compl. Analysis, V.2016.

\bibitem{8} D. F. Abzalilov, E. A. Shirokova. The approximate conformal mapping onto
simply and doubly connected domains. Complex Variables and Elliptic
Equations, V.62, pp. 554-565, (2017).

\bibitem{9} F. D. Gakhov. Boundary value problems. Pergamon Press, Oxford,
 1966.

\bibitem{10} A. W. K. Sangawi, A. H. M. Murid, M. M. S. Nasser. Annulus with circular slit map of bounded multiply connected regions via integral equation method. Bull. Malays. Math. Sci. Soc., V.35, pp. 945--959, (2012).

\bibitem{11} A. A. M. Yunus, A. H. M. Murid, M. M. S. Nasser,  Numerical conformal mapping and its inverse of unbounded multiply connected regions
onto logarithmic spiral slit regions and rectilinear slit regions, Proc.
of the Royal Society A -- Math. Phys. and Eng. Sci.,  vol. 470
(2162), Article No. 20130514, (2014).

\bibitem{12} A. A. M. Yunus, A. H. M. Murid, M. M. S. Nasser, Numerical evaluation of conformal mapping and its inverse for unbounded multiply
connected regions. Bull. Malays. Math. Sci. Soc.,  vol. 1 (24),
pp. 1--24, (2014).

\bibitem{13}R.  Wegmann, Methods for numerical conformal mappings. In Handbook of complex analysis: geometric function theory, vol. 2, Amsterdam, The Netherlands: Elsevier, pp. 351--477, (2005).
\bibitem{14} E. A. Shirokova, On the approximate conformal mapping of the unit
disk on a simply connected domain. Russian Math.,  vol. 58,
pp. 47--56, (2014).

\bibitem{15} E. A. Shirokova, P. N. Ivanshin, Approximate Conformal Mappings and Elasticity Theory. J. of Compl. Analysis, 2016 DOI:
10.1155/2016/4367205.

\bibitem{16} D. F. Abzalilov, E. A. Shirokova, The approximate conformal mapping
onto simply and doubly connected domains. Complex Variables and
Elliptic Equations, vol. 62, pp. 554--565, (2017).

\bibitem{BB} C.J. Bishop, Conformal mapping in linear time, Discrete and Comput. Geometry, vol 44, no. 2, pp. 330--428 (2010).

\bibitem{Leh} R. S. Lehman, Development of the mapping function at an analytic corner, Pacific J.
Math., 7, pp. 1437–1449 (1957).

\bibitem{Pm2} D.M. Hough, N. Papamichael,  The use of splines and singular functions in an integral equation method for conformal mapping. Numer. Math. 37, 133–147 (1981).

\bibitem{Iar} P. N. Ivanshin,
  Continued fractions and conformal mappings for domains with angle points, Arxiv:1711.04409


\bibitem{MZ}  A. Zygmund, Trigonometric Series, second ed., Cambridge University Press, London, New York, 1959.

\bibitem{J2}  D. Jackson, On the accuracy of trigonometric interpolation, Trans. Amer. Math. Soc. 14 (1913) 453–461.


\bibitem{17} T. K.  DeLillo, The Accuracy of Numerical Conformal Mapping Methods: A Survey of Examples and Results. SIAM J. Num. Anal.,
vol. 31, (1994).
\bibitem{18} T. K. DeLillo, On some relations among numerical conformal mapping
methods. J. Comput. Appl. Math., vol. 19, pp. 363--377, (1987).


%\bibitem{abz-shir}


\end{thebibliography}
\end{document}